\documentclass[12pt,a4paper]{amsart}
\usepackage{hyperref}
\usepackage{amssymb}
\usepackage{amsthm}
\usepackage{xypic}
\xyoption{all}

\newcommand{\A}{\ensuremath{\mathcal{A}}}

\newcommand{\D}{\ensuremath{\mathfrak{D}}}

\newcommand{\I}{\ensuremath{\mathcal{I}}}
\newcommand{\J}{\ensuremath{\mathcal{J}}}
\newcommand{\K}{\ensuremath{\mathcal{K}}}
\newcommand{\LLL}{\ensuremath{\mathcal{L}}}
\newcommand{\M}{\ensuremath{\mathcal{M}}}

\newcommand{\Q}{\ensuremath{\mathfrak{Q}}}
\newcommand{\R}{\ensuremath{\mathcal{R}}}
\newcommand{\SSS}{\ensuremath{\mathcal{S}}}
\newcommand{\T}{\ensuremath{\mathcal{T}}}

\newcommand{\V}{\ensuremath{\mathcal{V}}}
\newcommand{\W}{\ensuremath{\mathcal{W}}}
\newcommand{\X}{\ensuremath{\mathcal{X}}}

\newcommand{\CC}{\ensuremath{\mathbb{C}}}
\newcommand{\FF}{\ensuremath{\mathbb{F}}}

\newcommand{\PP}{\ensuremath{\mathbb{P}}}

\newcommand{\ZZ}{\ensuremath{\mathbb{Z}}}

\newcommand{\h}{\ensuremath{\mathfrak{h}}}
\newcommand{\s}{\ensuremath{\mathfrak{s}}}

\def\hol{{\mathcal{O}}}

\DeclareMathOperator{\Ext}{Ext}
\DeclareMathOperator{\Fix}{Fix}
\DeclareMathOperator{\Hom}{Hom}

\DeclareMathOperator{\Proj}{Proj}
\DeclareMathOperator{\Sym}{Sym}

\DeclareMathOperator{\coker}{coker}

\DeclareMathOperator{\supp}{supp}

\newtheorem{teo}{Theorem}[section]
\newtheorem{lem}[teo]{Lemma}

\newtheorem{cor}[teo]{Corollary}
\newtheorem{prop}[teo]{Proposition}

\theoremstyle{definition}
\newtheorem{df}[teo]{Definition}

\theoremstyle{remark}
\newtheorem{remark}[teo]{Remark}
\newtheorem{example}[teo]{Example}

\title{On surfaces with a canonical pencil}
\author{Roberto Pignatelli}
\address{Dipartimento di Matematica, Universit\`a di Trento. Via Sommarive 14,
loc. Povo, I-38050 Trento (Italy)}
\begin{document}

\begin{abstract}

We classify the minimal surfaces of general type with $K^2 \leq
4\chi-8$ whose canonical map is composed with a pencil, up to a
finite number of families.

More precisely we prove that there is exactly one irreducible family
for each value of $\chi \gg 0$, $4\chi-10 \leq K^2 \leq 4\chi-8$. All
these surfaces are complete intersections in a toric $4-$fold and bidouble 
covers of Hirzebruch surfaces. The surfaces with $K^2=4\chi-8$ were previously 
constructed by Catanese as bidouble covers of $\PP^1 \times \PP^1$.
\end{abstract}

%\thanks{The research of the authors was performed in the realm of the
%DFG SCHWERPUNKT "Globale Methode in der komplexen Geometrie", and of the
%EAGER EEC Project. The third author was  supported by the Schwerpunkt
%and by P.R.I.N. 2002 "Geometria delle variet\`a algebriche" of
%M.I.U.R. and is a member of G.N.S.A.G.A. of I.N.d.A.M. }

\maketitle
%\tableofcontents
\pagestyle{myheadings}
%\markboth{$K_S^2=4$, $p_g=q=1$}{$K_S^2=4$, $p_g=q=1$}
%\markright{On surfaces of general type with a canonical pencil and small $K^2$.}

\section*{Introduction}

The surfaces of general type form the class of the surfaces which have not yet been 
classified. One of the most successful approaches to study these 
surfaces has been to establish some inequality among the numerical 
invariants, and then to characterize the surfaces for which the inequality 
is an equality, or the difference of the two sides of the inequality is small.

The first famous inequality is the Noether inequality
$$K_S^2 \geq 2p_g(S)-4$$
which holds on any minimal surface of general type $S$.
As usual we denote by $K_S$ a canonical divisor, by $p_g(S)$ the dimension of the vector space $H^0(\hol_S(K_S))\cong H^0(\Omega^2_S)$, and we say that a surface $S$ is minimal if it is a minimal model, equivalently if $K_S$ is big and nef. 
We will later use also the invariants $q(S):=\dim H^1(\hol_S)$ and $\chi:=\chi(\hol_S)=1-q+p_g$. In the 70's the surfaces with $2p_g-4\leq K^2\leq 2p_g-2$ were object of a deep and successful investigation by Horikawa in the important papers \cite{horI}, \cite{horII}, \cite{horIII}, \cite{horIV}.

The second important example is the Castelnuovo inequality
 $$K_S^2 \geq 3p_g(S)-7$$
which holds on the minimal surfaces of general type whose canonical map, which is the rational map associated to $H^0(\Omega^2_S)$, is birational.
The surfaces for which the Castelnuovo inequality is an equality were investigated by Ashikaga and Konno in \cite{ak}.

This paper studies the surfaces of general type such that the image of the  canonical map is a curve. In this case one says that {\it the canonical map is
composed with a pencil}, since, by Stein factorisation, the
canonical map factors through a rational map with connected fibres onto a
smooth curve, which is said to be a {\it canonical pencil}. The {\it genus of the pencil} $g$ denotes by definition the genus of the general fibre.

In the celebrated paper \cite{bea}, Beauville studies the more general class of the surfaces whose canonical map is not birational.
Beauville's approach is ``up to a bounded family": he shows that if we
assume $\chi$ big enough, which is
equivalent to disregarding a finite number of components of the moduli
space of the surfaces of general type, the situation simplifies
considerably.

In our case \cite{bea} proves
\begin{teo}[{\bf Beauville}]\label{beauville}
If $S$ is a minimal surface of general type whose canonical map is composed
with a pencil and $\chi(\hol_S)>20$ then the general fibre has genus $2\leq g \leq 5$.

Moreover the pencil is free, {\it i.e.} the canonical map is a morphism. In other words, writing the canonical system as $|K_S|=\Phi+|M|$, the sum of a fixed part and a movable part, then the movable part $|M|$ has no base points.
\end{teo}

Beauville first proves that the inequality $K^2 \geq 2(g-1)(\chi-2)$ holds,
and then the theorem follows by the Miyaoka--Yau inequality $K^2 \leq
9\chi$.
Earlier Reid had asked in \cite[Problem R.5]{copenaghen} whether in the case $g=2$ the inequality $K^2 \geq 4p_g-6$ holds. 
This was proved by \cite{xiao} as follows

\begin{teo}[{\bf Xiao}]\label{Xiao}
If $S$ is a minimal surface of general type whose canonical map is
composed with a free pencil of curves of genus $2$ then $$K^2 \geq
4p_g-6 \geq 4\chi -10.$$

The inequality is sharp, since there are examples of surfaces with $q=0$
as above for all values of the pair of invariants $(K^2,p_g)$ in the
range $p_g \geq 2$, $4p_g-6 \leq K^2 \leq 4p_g$.
\end{teo}

We give a proof of Xiao's inequality $K^2 \geq
4p_g-6$ in Theorem~\ref{gang}.

A natural question arises: is it possible to classify the minimal surfaces 
of general type whose canonical map is composed with a pencil and with $K^2 - (4\chi
-10)$ ``small enough"? To the best of our knowledge, nobody has
studied this problem up to now.

In this note we answer this question. Our result is
the following

\begin{teo}\label{g=2}
Let $S$ be a minimal surface of general type with $K_S^2 \leq
4\chi-8$. If the canonical system is a free pencil of genus $2$ curves, then $S$
is regular and the pencil has rational base.

If we further assume either $K^2 \leq 4\chi-10$ or $K^2 <5(\chi-3)$
then moreover $S$ is birational to a complete intersection $X=Q\cap G$
in the toric $4-$fold $\PP:= \Proj(\Sym V)$ where
$$
V:=\hol_{\PP^1}(1)x_0 \oplus \hol_{\PP^1}(\chi)x_1 \oplus
\hol_{\PP^1}(K^2-2\chi+8)y \oplus \hol_{\PP^1}(K^2-\chi+7)z
$$
where the grading of $\Sym V$ is given by $\deg x_i=1,\ \deg y=2,\ \deg
z=3$.

$X$ has only canonical singularities.  $Q,G$
are of the form
$$
Q:=\left\{ x_0^2+q_x x_1^2+q_y y=0\right\},
$$
$$
G:=\left\{ z^2 +\sum_{\substack{i,j,k \geq 0\\
i+j+2k=6}} G_{ijk}x_0^ix_1^jy^k=0\right\},
$$
where $q_x, q_y$ and $G_{ijk}$ are homogeneous polynomials on $\PP^1$.

Conversely, for any pair of positive integers $(K^2, \chi)$, the minimal resolution $S$ of the singularities of a complete intersection $X$ as
above is a minimal surface of general type with $K_S^2=K^2$,
$\chi(\hol_S)=\chi$, $q(S)=0$ and the canonical map of $S$ is the composition
$ S \stackrel{f}{\rightarrow}\PP^1
\stackrel{Ver}{\hookrightarrow}\PP^{\chi-2}$ where $f$ is the genus $2$ pencil
induced by the projection $\PP \rightarrow \PP^1$ and $Ver$ is the
$(\chi-2)^{th}-$Veronese embedding.

In particular, if $(K^2,\chi)$ are integers with $\chi >>0$ and
$4\chi-10 \leq K^2 \leq 4\chi-8$, then the closed subset of the moduli
space of the minimal surfaces of general type given by the surfaces
$S$ with $K^2_S=K^2$, $\chi(\hol_S)=\chi$ whose canonical map is
composed with a pencil is irreducible of dimension $12\chi-2K^2-15$.

$X$ is the canonical model of $S$ unless $K^2_S=2$, $\chi(\hol_S)=3$.
\end{teo}

We assume that the pencil is free: this is not a very restrictive assumption. Indeed it is well known (see \cite{horP} or \cite[Proposition 4.1]{xiao} that the only non free pencil of genus $2$ curves on a surface of general type is the canonical system  of a hypersurface of degree $10$ in $\PP(1:1:2:5)$.

We do not completely classify the case $K^2\leq 4\chi-8$ because we have further assumed $K^2\leq 4\chi-10$ or $K^2< 5(\chi-3)$. Also this assumption is not very restrictive: we are excluding only the cases $K^2=4\chi-9$, $\chi \leq 6$ and $K^2=4\chi-8$, $\chi \leq 7$.  
The assumption $K^2< 5(\chi-3)$ is necessary: in Example~\ref{example} we construct  surfaces with $(K^2,\chi)=(15,6)$ and $(20,7)$ whose canonical map is composed with a pencil of genus $2$ curves which are not considered in Theorem~\ref{g=2}. For these small values of the invariants there will be some exceptional families. We think that the strategy of Example~\ref{example} can lead to a complete classification of the case $K^2\leq 4\chi-8$.

The construction in Theorem~\ref{g=2} gives also surfaces with
canonical map composed with a pencil and $K^2 > 4\chi-8$; at least for
every pair $(K^2,\chi)$ with $\chi \geq 3$, $4\chi-10 \leq K^2 \leq
4\chi-6$. Anyway, when $K^2 > 4\chi-8$ we can't guarantee that these
are all the surfaces with a canonical pencil. Indeed, using an
alternative description of these surfaces as bidouble covers, we can show
that if $K^2=4\chi-6$ then this family is always properly contained in
a larger irreducible family of surfaces with a canonical pencil.

The case $K^2_S=2$, $\chi(\hol_S)=3$ has been fully described in \cite[Theorem~1.3]{horIV}; 
in this case $X$ is not the canonical model since, if $\Phi$ is the fixed part 
of the canonical system, $K\Phi=0$.

If $g \geq 3$ Beauville's inequality implies $K^2 \geq 4\chi-8$ and 
even stronger inequalities have been proved later by Sun (see \cite{sun1}, \cite{sun2}; these results are summarized in Theorem~\ref{sun}; see also \cite{zuc}, \cite{YM}).
From them it follows that if $K^2$ is sufficiently near to $4\chi -10$ and $\chi$ is
sufficiently big, then the genus of the canonical pencil has to be
$2$. We can then prove

\begin{teo}\label{allgenera?}
Let $S$ be a minimal surface of general type with $K_S^2 \leq
4\chi(\hol_S)-8$ and canonical map composed with a pencil. Assume moreover
$\chi(\hol_S) \gg 0$. Then $S$ is one of the surfaces in Theorem~\ref{g=2}. 
\end{teo}

\

We discuss the strategy of the proof.

By the above mentioned results we can assume that the
pencil $f\colon S \rightarrow B$ is free, $B$ is rational and the fibres have genus $g=2$. We
classify these, for $\chi \gg 0$ and $K^2 \leq 4\chi -8$, describing
their relative canonical algebra 
$$\R(f):= \bigoplus_n \R_n:= f_* (\omega_S \otimes f^*\omega_B^{-1})^{\otimes n}.$$

The main difficulty lies in the fact that these genus $2$ fibrations have big
Horikawa number $H:=K^2-2\chi+6$ (cf. \cite{horP}, \cite{rei}): it
grows asymptotically as
$2\chi$. When studying the relative canonical algebra $\R$ of genus $2$
fibrations over surfaces with fixed invariants $(K^2,\chi)$, big values
of $H$ correspond to many {\it a priori} possibilities for its second
graded piece $\R_2$, which is the key bundle to be computed.
Despite that, we prove that, if the
fibration is canonical, $K^2 \leq 4\chi -8$ and $\chi \gg 1$, then
$\R_2$ is determined up to isomorphism, while this is not true for
small $\chi$, when some {\it exceptional} families appear.

We explain the geometric idea underlying our proof.
Any surface $S$ with a genus $2$ fibrations is birational to the relative 
minimal model $X$ of the fibration. $X$ has an involution $i$, 
given by the hyperelliptic involutions of the fibres, such that the quotient $C=X/i$ 
has a natural structure of conic bundle. The horizontal part $\h$ of the
fixed part of $|K_X|$ cuts a canonical divisor on the general fibre, 
and therefore it is invariant and maps to a section $\s$ of the conic bundle $C$.
This section determines a subsheaf $\SSS$ of the relative canonical algebra.
We show that $\s$ can't be contained in the branch divisor $\Delta$ of the double 
cover $X\rightarrow C$ by studying $\SSS$ and the isolated branch points of the same map. 
This implies the nontriviality of certain maps between line bundles on the base curve. 
A first consequence is a short proof of Xiao's inequality.

The key idea is to consider the first infinitesimal neighbourhood of $\s$, {\it i.e.} 
the subscheme $2\s \subset C$, which determines a new sheaf of
algebras $\SSS'$. Interpreting the condition $\s \not\subset \Delta$ in terms of $\SSS'$ 
we were able to compute first $\R_2$, and then the whole sheaf of algebras $\R$.

\

The paper is structured as follows.

In Section~\ref{families} we study the surfaces constructed in
Theorem~\ref{g=2} as complete intersections in a toric
$4-$fold, showing that the family is not empty for every
value of the pair $(\chi, K^2)$ in the range $\chi \geq 3$, $4\chi-10
\leq K^2 \leq 4\chi-6$. Then we check that these surfaces have
the properties stated in Theorem~\ref{g=2}. We also show 
how to deform these surfaces to surfaces which still have the genus $2$ pencil
but whose canonical map is not composed with it, so showing that none of our 
families is an irreducible component of the moduli space of surfaces of general type.
Moreover we note that, if $\chi$ is big enough, these surfaces are bidouble covers of 
a Hirzebruch surface $\FF_{|4\chi-8-K^2|}$.

In Section~\ref{preliminaries} we recall some results from the theory of genus
$2$ fibrations, and we introduce the sheaf of algebras $\SSS$ associated to a section
$\s$ of the related conic bundle $C \rightarrow B$.

In Section~\ref{s} we study $\SSS$, and we give a short proof of
Xiao's inequality. 

In Section~\ref{2s} we introduce the sheaf of algebras
$\SSS'$ of $2\s$ and use it (and a lifting lemma) to prove Theorem~\ref{main}
which determines exactly for which pair $(K^2,\chi)$ there are surfaces with
canonical map composed with a free genus $2$ pencil not in the families
described in Theorem~\ref{g=2}. Theorem~\ref{g=2} and Theorem~\ref{allgenera?} follow then
easily.

{\bf Acknowledgements.}
I'm indebted to Miles Reid for suggesting the problem to me and for
many enlightening discussions when I was guest of the University of
Warwick in the occasion of the WAG 07-08 Symposium, where most of the
results of this paper were obtained. To him belongs the idea of 
Remark~\ref{miles}.
I thank the referee for his useful comments and suggestions.
%\section*{Notation}

\section{Some surfaces with a canonical pencil}\label{families}

In this section we study the surfaces constructed in Theorem~\ref{g=2} as complete intersection in $\PP:=\Proj(\Sym V)$ where $V$ is a vector bundle depending on two integers $\chi$ and $K^2$. For later convenience we introduce the integers $p_g:=\chi-1$ and $\Theta:=K^2-4\chi+10$. We will prove later in this section that $\chi$, $p_g$ and $K^2$  equal resp. $\chi(\hol_S)$, $p_g(S)$ and $K^2_S$ (in particular $q(S)=0$).

Then
$$
V=\hol_{\PP^1}(1)x_0 \oplus \hol_{\PP^1}(p_g+1)x_1 \oplus \hol_{\PP^1}(2p_g+\Theta)y \oplus \hol_{\PP^1}(3p_g+\Theta)z.
$$

$\PP$ is the $4-$fold $\Proj(\Sym V)$ where the grading of $\Sym V$, sheaf of graded algebras
over $\PP^1$, is given by $\deg x_i=1,\ \deg y=2,\ \deg z=3$.
So the first graded pieces of $\Sym V$ are
$(\Sym V)_1=\hol_{\PP^1}(1)x_0 \oplus \hol_{\PP^1}(p_g+1)x_1$,
$(\Sym V)_2=\hol_{\PP^1}(2)x_0^2 \oplus \hol_{\PP^1}(p_g+2)x_0x_1
\oplus \hol_{\PP^1}(2p_g+2)x_1^2\oplus \hol_{\PP^1}(2p_g+\Theta)y$,
and so on.

$\PP$ is a bundle over
$\PP^1$ whose fibres are weighted projective spaces $\PP(1:1:2:3)$. In
$\PP$ we have the hypersurfaces
$
Q=\left\{ x_0^2+q_x x_1^2+q_y y=0\right\}$, $G=\left\{ z^2 +\sum_{i+j+2k=6} G_{ijk}x_0^ix_1^jy^k=0\right\}
$
where the $q_i$ and the $G_{ijk}$ are homogeneous polynomials
on $\PP^1$ of the right degree to have $Q$ and $G$ well defined:
$$\begin{aligned}
\deg q_x&=2p_g\\
\deg q_y&=2p_g-2+\Theta\\
\deg G_{ijk}&=[i+(p_g+1)j+(2p_g+\Theta)k]-2(3p_g+\Theta)\\
&=(j+2k-6)p_g+(k-2)\Theta+(i+j)\\
&=-ip_g+(k-2)\Theta+(6-2k)
\end{aligned}$$

$X:=Q \cap G$ is a surface and the
fibres of the surjective morphism $f\colon X \rightarrow \PP^1$
are complete intersections of type $(2,6)$ in $\PP(1:1:2:3)$.
We assume that $X$ has only canonical singularities, which are isolated singularities, and therefore the general fibre of $f$ is smooth. Then $f$ is a genus $2$ fibration.

\

\paragraph{\it The canonical system of $S$, the relative canonical algebra of $f$ and the invariants.}

By adjunction, denoting by $H$ the class of $|\hol_{\PP}(1)|$, the
linear system of the relative hyperplanes
defined by maps $\hol_{\PP^1} \hookrightarrow (\Sym V)_1$, and by $F$ the
pull back of a point of $\PP^1$
\begin{align*}
K_{X|\PP^1}&=(K_{\PP|\PP^1}+Q+G)_{|X}\\
           &=((6p_g+2\Theta+2)F-7H-2F+2H-(6p_g+2\Theta)F+6H)_{|X}\\
&=H_{|X}.
\end{align*}

Consider the relative canonical algebra
$$
\R(f):=\bigoplus_{d\geq 0} f_{*} \hol_X(dK_{X|\PP^1}).
$$

Since $K_{X|\PP^1}=H_{|X}$, $\R(f)$ is the quotient of $\Sym V$ by the sheaf of ideals $\I$ generated by the polynomials defining $Q$ and $G$.

In particular $K_X=(H-2F)_{|X}$, the restriction map
$H^0(\PP, H-2F)\rightarrow  H^0(X, K_X)$ is an isomorphism and
$H^0(X, K_X)$ is cut by the hypersurfaces $\{hx_1=0\}$,
with $h$ varying in $H^0(\hol_{\PP^1}(p_g-1))$: the fixed part of the canonical system is cut by $\{x_1=0\}$ and the canonical map is the
composition of $f$ with the $(p_g-1)-$uple embedding $\PP^1
\rightarrow \PP^{p_g-1}$.

Then $p_g(S)=p_g$, and by the standard formulae for genus 2 fibrations (see \cite[Theorem~ 4.13]{cp}), $\chi(\hol_S)=\deg f_*(K_{X|\PP^1})-1=1+p_g+1-1=\chi$. Moreover
\begin{align*}
K^2_S&=(H-2F)^2\cdot Q \cdot G\\
&=12H^4-4(3p_g+\Theta+15)H^3F\\
&=2\left( 1+p_g+1+\frac{2p_g+\Theta}2+\frac{3p_g+\Theta}3\right) -4\frac{3p_g+\Theta+15}6\\
&=4p_g-6+\Theta=K^2.
\end{align*}

\

\paragraph{\it The family.}

\begin{prop}\label{notempty}
If $0 \leq \Theta \leq 4$ then the general $X$ is smooth.
\end{prop}

\begin{Proof}
The equation of $Q$ does not involve the variable $z$, so $Q$
is a cone over $C:=Q \cap \PP'$, where $\PP':=\{z=0\}$ is a
$\PP(1:1:2)-$bundle over $\PP^1$. Since $p_g\geq 2$ and $\Theta \geq
0$, $q_x$ and $q_y$ have positive degree and therefore by Bertini's theorem
the general $C$ is {\it quasi-smooth} and
more precisely has only singularities of type $A_1$ (nodes) at the
intersection with the section $\{x_0=x_1=0\}$: $2p_g-2+\Theta$ nodes
dominating the zeroes of $q_y$.

By the equation of $G$, $X$ is a double cover of $C$ branched at
these nodes and along the curve $\Delta:=C \cap \bar{\Delta}$ where
$\bar{\Delta}:=\{ \sum  G_{ijk}x_0^ix_1^jy^k=0\}$ is a surface in $\PP'$.
To prove the smoothness of the general $X$ we need that for general
choice of the $G_{ijk}$, $\Delta$ is smooth and does not contain the nodes.

We study the fixed locus of $|\Delta|$. Note that $|\bar{\Delta}|$  has a big fixed locus: in fact if $p_g \gg 0$
all $G_{ijk}$ with $i>0$ vanish (having negative degree) and therefore the curve $\{x_1=y=0\}$ is contained in $\Fix(|\bar{\Delta}|)$.

If $0 \leq \Theta \leq 3$, since $\deg G_{003}=\Theta$ and $\deg
G_{060}=6-2\Theta$, the coefficients of both the monomials $y^3$ and
$x_1^6$ can be chosen not identically zero and therefore $\Fix(|\bar{\Delta}|$) is contained in the curve $\{x_1=y=0\}$. Since
this curve does not intersect $C$ (due to the term $x_0^2$ in the
equation of $Q$), $|\Delta|$ has no base points and its general element is
smooth and irreducible.

If $\Theta=4$, $G_{060}=0$ and therefore, if $p_g \gg 0$, the equation of $G$ is divisible by $y$: $\{y=0\}\cap C$ is a smooth fixed component of $|\Delta|$. Since $\deg G_{041}=0$, the residual linear system has no base points (as in the previous case)
and its general element does not intersect the fixed component, so $\Delta$ is
disconnected but still smooth.

\qed
\end{Proof}

\begin{prop}
If $0 \leq \Theta \leq 2$ and $p_g>6-2\Theta$  the subset of the moduli space of surfaces of general type given by these surfaces is not empty, unirational of dimension $4p_g+9-2\Theta$.
\end{prop}

\begin{Proof}The family is not empty by Proposition~\ref{notempty}. The unirationality follows since the family is parameterized by the coefficients of  general polynomials on $\PP^1$ of given degrees.
The dimension count is identical to the analogous one in
\cite[Proposition~5.2]{mealone}.

\qed
\end{Proof}

One can use the same method to compute the dimension of the other families, which can be a little bit bigger. Indeed, if $\Theta \leq 2$ and $p_g > 6-2\Theta$, then no term including the variable $x_0$ appears in the equation of $G$: when these terms appear, they give a slightly larger number of moduli.

The forthcoming Remark~\ref{miles} shows that none of these families is an irreducible component 
of the moduli space. This was already known for $p_g=2$, $K^2=2,3$ by \cite{horIV}, resp. \cite{fo}. 
In general, we do not know the codimension of these families in the irreducible component of the moduli space containing them.

\begin{remark}\label{miles}
Reid showed me how to deform all the surfaces in Theorem~\ref{g=2} to surfaces which still have the genus $2$ pencil but whose canonical map is no longer composed with it.

Consider
$$
V'=\hol_{\PP^1}(1)x'_0 \oplus \hol_{\PP^1}(2)x'_1 \oplus \hol_{\PP^1}(p_g)x'_2
\oplus \hol_{\PP^1}(2p_g+\Theta)y' \oplus \hol_{\PP^1}(3p_g+\Theta)z'.
$$
and $\PP':=\Proj(\Sym V')$ where the grading of $\Sym V'$ is given by $\deg x'_i=1,\ \deg y'=2,\ \deg z'=3$. In $\PP' \times \CC$ we take three hypersurfaces
$$
L=\left\{ \lambda x'_0-t_0x_1'+t_1^{p_g-1}x_2'=0\right\},$$
$$
Q'=\left\{ (x'_0)^2+q'_{11} (x'_1)^2+q'_{12} x'_1x'_2+q'_{22} (x'_2)^2+q'_y y'=0\right\},$$ $$G'=\left\{ (z')^2 +\sum_{i+j+k+2l=6} G'_{ijkl}(x'_0)^i(x'_1)^j(x_2')^k(y')^l=0\right\}$$
where $\lambda$ is a coordinate on $\CC$, $(t_0:t_1)$ are homogeneous coordinates on $\PP^1$, $q'_y$, the $q'_{ij}$ and the $G'_{ijkl}$ are homogeneous polynomials on $\PP^1$ of the right degree to have $L$, $Q'$ and $G'$ well defined.

Consider then the 3-fold $\X:= L \cap Q' \cap G' \subset \PP' \times \CC$ and the map $\xi \colon \X \rightarrow \CC$ restriction of the projection on the second factor $\xi_2 \colon \PP' \times \CC \rightarrow \CC$.

Note that $\xi_2^{-1}(0) \cap L \cong \PP$, the isomorphism being given 
by $x_0'=x_0$, $x_1'=t_1^{p_g-1}x_1$, $x'_2=t_0x_1$,  $y=y'$, $z=z'$. 
Then, setting $Q:=Q' \cap \xi_2^{-1}(0) \cap L \subset \PP$ and $G:=G' \cap \xi_2^{-1}(0) \cap L \subset \PP$ we see that
 $X_0:=\xi^{-1}(0)=Q \cap G$ is one of the surfaces in Theorem~\ref{g=2}. The reader can easily 
check that every surface in Theorem~\ref{g=2} can be obtained in this way.

Consider now a general element of the family $X_{\lambda}:=\xi^{-1}(\lambda)$. The natural projection onto $\PP^1$ 
is still a free pencil of genus $2$ curves. A basis for $H^0(K_{X_{\lambda}}) \cong \CC^{p_g}$ is $\{x_1', t_0^{p_g-2}x'_2,
t_0^{p_g-3}t_1x'_2, \ldots , t_1^{p_g-2}x'_2\}$. 
In particular the base locus of $|K_{X_{\lambda}}|$ is $X_{\lambda} \cap\{x_1'=x_2'=0\}$.
If $\lambda \neq 0$ from the equation of $L$ follows $x_0'=0$ and from the equations of $G'$ 
we deduce that the base locus of $|K_{X_{\lambda}}|$ is finite. In particular the canonical map is not composed with any free pencil.
\end{remark}

\

\paragraph{\it Description as bidouble covers}

For $p_g>6$, $x_0$ does not appear
in the equation of $G$. In this case $X$ is invariant by the $(\ZZ/2\ZZ)^2$ action on $\PP$ defined by
$
(x_0:x_1:y:z) \mapsto (\pm x_0:x_1:y:\pm z).
$

The quotient
$X/(\ZZ/2\ZZ)^2$ is the $\PP(1:2)$ bundle over $\PP^1$ given
by the variables $x_1$ and  $y$, which is the Hirzebruch surface
$$\FF_{|2-\Theta|}=\Proj(\Sym(\hol(2p_g+2)x_1^2\oplus
\hol(2p_g+\Theta)y)).$$

As in \cite{cat99}, a bidouble cover of a surface is
determined by three divisors $D_1$, $D_2$, $D_3$ such that the branch divisors   of the $3$ intermediate double covers are $D_1+D_2$, $D_1+D_3$ and $D_2+D_3$. In our case $D_1=\{q_xx_1^2+q_yy=0\}$,
$D_2=\{\sum G_{0jk}x_1^jy^k=0\}$, $D_3=\{x_1=0\}$.

Writing
$\Gamma_{\infty}$ for the only negative section of the ruling
$|\Gamma|$ of $\FF_r$, and $|\Gamma_1|$, $|\Gamma_2|$ for the two
rulings of $\PP^1 \times \PP^1$, the $D_i$ are general in the following linear system

\

\renewcommand{\arraystretch}{1.3}
\begin{tabular}{|c|c|c|c|c|c|}
\hline
$\Theta$&base&$|D_1|$&$|D_2|$&$|D_3|$\\
\hline
\hline
$0$&$\FF_2$&$|\Gamma_{\infty}+2p_g\Gamma|$&$|3\Gamma_{\infty}+6\Gamma|$&$\Gamma_{\infty}$\\
$1$&$\FF_1$&$|\Gamma_{\infty}+2p_g\Gamma|$&$|3\Gamma_{\infty}+4\Gamma|$&$\Gamma_{\infty}$\\
$2$&$\PP^1 \times \PP^1$&$|\Gamma_1+2p_g\Gamma_2|$&$|3\Gamma_1+2\Gamma_2|$&$|\Gamma_1|$\\
$3$&$\FF_1$&$|\Gamma_{\infty}+(2p_g+1)\Gamma|$&$|3\Gamma_{\infty}+3\Gamma|$&$|\Gamma_{\infty}+\Gamma|$\\
$4$&$\FF_2$&$|\Gamma_{\infty}+(2p_g+2)\Gamma|$&$\Gamma_{\infty}+|2\Gamma_{\infty}+4\Gamma|$&$|\Gamma_{\infty}+2\Gamma|$\\
\hline
\end{tabular}

\

This table shows in particular that the surfaces we have constructed with
$\Theta=2$ are the same as those constructed by Catanese in \cite[Example~2]{cat99}, (the case $a=2$, $k=0$ in the notation there). Catanese in
the same Example gave also surfaces with $\Theta=4$ or $6$ as bidouble
covers of $\PP^1 \times \PP^1$: it is easy to see that our examples
with $\Theta=4$ are a degeneration of Catanese's examples, obtained
by deforming $\PP^1 \times \PP^1$ to $\FF_2$.

\begin{prop}
There are surfaces whose canonical map is composed with a genus $2$
pencil which are bidouble covers of a Hirzebruch surface for each
value of $p_g \geq 2$, $\Theta:=K^2-4p_g+6 \in \{0,1,2,3,4,5,6\}$.
\end{prop}

\begin{Proof}
The above table proves the statement for $\Theta \leq 4$, and
Catanese's examples above mentioned the case $\Theta=6$. For
$\Theta=5$ it is enough to take a bidouble cover of $\FF_1$ with
branch divisors $D_1 \in |\Gamma_{\infty}+(2p_g+2)\Gamma|$, $D_2 \in
\Gamma_{\infty}+ |2(\Gamma_{\infty}+\Gamma)|$, $D_3 \in
|\Gamma_{\infty}+2\Gamma|$.

\qed
\end{Proof}

The existence of surfaces with these values of the invariants was already known
to Xiao who constructed in \cite{xiao} one example for each value of
$p_g$, $K^2$ in the above range.

\section{Genus 2 fibrations, related conic bundles and their
  sections}\label{preliminaries}

If a surface $S$ has a fibration onto a curve $f\colon S \rightarrow B$ with fibres of genus $g \geq 2$, then $S$ is birational to
$X:=\Proj(\R(f))$, where $\R(f)$ (or $\R$ for short) is the relative canonical algebra
$$
\R:=\bigoplus_{d\geq 0} f_{*} \hol_S(dK_{S|B}).
$$

$X$, the {\it relative canonical model} of $f$, is obtained from $S$
by contracting all {\it vertical} $(-2)-$curves. To simplify the
arguments, it is convenient to replace $S$ by $X$, considering (with a
mild abuse of notation) the fibration $f \colon X \rightarrow B$.

The hyperelliptic involution of the fibres extends to an involution $i$ of $X$, inducing a quotient surface $C:=X/i$. The fibration $f$ factors as $\pi \circ \gamma$, where $\pi \colon C \rightarrow B$ is a conic bundle and $\gamma \colon X \rightarrow C$ is a finite double cover branched at a finite set of points $P \subset C$ and along a curve $\Delta \subset C$. Obviously
\begin{equation}\label{PDelta}
P \cap \Delta =\emptyset.
\end{equation}

 Since $X=\Proj(\R)$, then $C=\Proj(\A)$ for a subalgebra $\A \subset \R$ which has been studied in  \cite{cp} (see also \cite{rei}).

\begin{lem}[\cite{cp}]\label{A}
$\A$, as sheaf of algebras, is generated by $\A_1=\R_1$ and $\A_2=\R_2$ and
related by the multiplication map
$$
\sigma_2 \colon \Sym^2 \A_1 \hookrightarrow \A_2
$$
whose cokernel is the structure sheaf of an effective divisor $\tau$
on $B$ of degree $K_X^2-2\chi+6$. Moreover $\deg \R_1=\chi(\hol_S)-b+1$ where $b$ is the genus of $B$.
\end{lem}

\begin{Proof}
This was proved in \cite{cp}. More precisely generators and
relations of $\A$ are computed by Lemma~4.4 and Remark~4.5 and the
cokernel of $\sigma_2$ is computed in Lemma~4.1.
For the degrees of $\tau$ and $\R_1$ see the formulae in Theorem~4.13.

\qed
\end{Proof}

In the next remark we recall the description of the stalks of $\A$ and $\R$ at
any point of the base curve $B$ which will be used in the
forthcoming sections.

\begin{remark}\label{tau}
If $p \in B$ does not belong to $\supp \tau$, then the stalk of $\R$
at $p$ is isomorphic to
$$\hol_{B,p}[x_0,x_1,z]/(z^2-f_6(x_0,x_1;t) );$$
else, letting $r$ be the multiplicity of $p$ in $\supp \tau$, the
stalk of $\R$ at $p$ is isomorphic to
\begin{equation}\label{stalk}
\hol_{B,p}[x_0,x_1,y,z]/(t^ry-f_2(x_0,x_1;t),z^2-f_6(x_0,x_1,y;t) ).
\end{equation}

Here $f_d$ are weighted homogeneous polynomials of degree $d$, where the variables $x_i$ have weight $1$ and $y$ has weight $2$. The involution acts as $z \mapsto -z$ ($z$ has weight $3$) fixing all the other variables. $\A$ is the subalgebra generated by $x_0$, $x_1$ and, in the second case, $y$.

From this description follows that $\pi$ maps $P$ bijectively to $\supp \tau$
and, equivalently $f$ maps bijectively the set of the isolated fixed points of
$i$ to $\supp \tau$; in fact, if $p \in
\supp \tau$  there is an isolated fixed point contained in $f^{-1}(p)$ which is, in the local
coordinates induced by the above description of the stalk of $\R$ at $p$, the
point
$Q_0:=((x_0:x_1:y:z);t)=((0:0:1:0);0)$.

If $r=1$ and $f_2$ and $f_6$ are general, then $f^{-1}(p)$ has two irreducible components, two elliptic curves intersecting transversally in $Q_0$. Each of these curves is mapped by $\gamma$ to a rational curve: the two rational curves intersect in $P_0=\gamma(Q_0) \in P$, which is a node of $C$.
\end{remark}

By Remark~\ref{tau} follows, if we denote by $\R_3'$ the direct summand of
$\R_3$ of rank $1$ locally generated by $z$,
\begin{prop}[\cite{cp}]\label{Delta}
$\R \cong \A \oplus (\A \otimes \R_3') [-3]$ as graded $\A-$module.
The ring structure of $\R$ gives a multiplication map $\delta \colon (\R_3')^2 \rightarrow \A_6$
inducing the divisor $\Delta \subset C$.
Moreover $\R_3' \cong \det \A_1 \otimes \hol_B(\tau)$.
\end{prop}

\begin{Proof}
See \cite{cp}, in particular Proposition~4.8. Note that $\R_3'$ is
$V_3^+$ in the notation of \cite{cp}.

\qed \end{Proof}

\begin{df}
Let $f\colon S \rightarrow B$ be a genus $2$ fibration, and consider the
associated conic bundle $\pi \colon C=\Proj(\A) \rightarrow B$.

To any section $\eta$ of $\pi$ we associate a sheaf of graded algebras $\SSS$,
which is the quotient of $\A$ by the ideal sheaf of the elements vanishing
along the curve $s:=\eta(B)$.
\end{df}

Note that by definition each homogeneous piece $\SSS_d$ of $\SSS$ is of rank
$1$ and torsion free, so a line bundle.

\begin{prop}\label{s->tau'}
There exists an effective divisor $\tau' < \tau$ such that $\forall
d \geq 1$, $\SSS_d \cong \SSS_1^d\left( \left\lfloor \frac{d}2\right\rfloor
\tau' \right)$.
\end{prop}
\begin{Proof}
The statement is empty for $d=1$, so we start with $d=2$. By Lemma~\ref{A}
we have a commutative diagram with exact rows and columns
$$
\xymatrix{%
0\ar[r]&\Sym^2 \A_1\ar[r]\ar[d]&\A_2\ar[r]\ar[d]&\hol_\tau\ar[r] & 0 \\
0\ar[r]&\SSS_1^2\ar[r]\ar[d]&\SSS_2\ar[d]& &  \\
&0&0& &  \\
}
$$
It follows that there exists a surjection from
$\hol_{\tau}$ to the cokernel of the map $\SSS_1^2 \hookrightarrow
\SSS_2$, which is then isomorphic to $\hol_{\tau'}$ for some
$\tau'<\tau$. In particular $\SSS_2 \cong \SSS_1^2(\tau')$.

To conclude we show that, $\forall k \geq 2$,  $\SSS_{2k} \cong
\SSS_2^k$ and $\SSS_{2k-1}\cong \SSS_2^{k-1} \otimes \SSS_1$.

The composition of the maps $\Sym^k \A_2 \rightarrow
\A_{2k} \rightarrow \SSS_{2k}$ is surjective (since both maps
 are surjective) and factors through the multiplication map
$\SSS_2^k \rightarrow \SSS_{2k}$, which is therefore surjective too;
the injectivity follows since it is a map between line bundles. A similar
argument works for the map $\SSS_2^{k-1} \otimes \SSS_1 \rightarrow
\SSS_{2k-1}$.

\qed \end{Proof}

Recall that by Remark~\ref{tau}, $\pi_{|P}$ maps $P$ bijectively onto
$\supp \tau$. In the next remark  we show that, roughly speaking,
$\tau' < \tau$ detects those points of $P$ which belong to $s$.

\begin{remark}\label{tau'}
Let $p$ be a point of $\tau$ with multiplicity $r>0$; we write the
stalks of $\R$ and $\A$ at $p$  as in Remark~\ref{tau}, by asking
further that $x_1$ generates the kernel of the map $\A_1 \rightarrow \SSS_1$.

$\SSS_2$ is a line bundle and the surjection $\A_2 \rightarrow \SSS_2$ factors through $(\A/x_1)_2$, which has rank $1$ but possibly torsion. More precisely, if $p \in \supp \tau$ then the stalk $(\A/x_1)_{2,p}\cong (\hol_{B,p}[x_0,y]/(t^ry-f_2(x_0,0;t))_2$ has torsion when $t$ divides $f_2(x_0,0;t)$. It
follows that
$$\SSS_{2,p}\cong (\hol_{B,p}[x_0,y]/((t^ry-f_2(x_0,0;t))/t^{r''})_2$$
where $r''$ is the maximal power of $t$ which divides $t^ry-f_2(x_0,0;t)$.

The stalk at $p$ of $\SSS_1$ is generated by the class of
$x_0$. Therefore
$$
\hol_{\tau'} \cong \coker (\SSS_1^2 \rightarrow \SSS_2) \cong \bigoplus_{p \in \supp \tau} (\hol_{B,p}[y]/(t^{r-r''}y))_2;
$$
the multiplicity of $p$ in $\tau'$ is $r':=r-r''$.

Note that $p$ belongs to $\tau'$ if and only if $r \neq r''$ {\it
  i.e.} the section passes through $P_0$. Therefore, if $\tau$
is reduced, then $\tau'<\tau$ detects the points of $P$ contained in
$s$. Moreover
\begin{equation}\label{tau'not0}
\tau'=0 \Leftrightarrow P \cap s = \emptyset.
\end{equation}
\end{remark}

\section{The section of the conic bundle induced by Fix(K)}\label{s}

We want to prove the following
\begin{teo}\label{main}
Let $S$ be a minimal surface of general type with $K_S^2 \leq
4\chi-8$, having a free canonical pencil of genus $2$ curves. Then $S$
is regular and the pencil has rational base.

If moreover we assume either $\Theta:=K^2-4\chi+10=0$ or $p_g >\Theta+4$
then $\R \cong (\Sym V)/\I$, with
$$
V:=\hol_{\PP^1}(1)x_0 \oplus \hol_{\PP^1}(\chi)x_1 \oplus \hol_{\PP^1}(K^2-2\chi+8)y \oplus \hol_{\PP^1}(K^2-\chi+7)z,
$$
the grading of $\Sym V$ is given by $\deg x_i=1,\ \deg y=2,\ \deg
z=3$, and $\I$ is the sheaf of ideals generated by
$x_0^2+q_x x_1^2+q_y y, z^2 +\sum_{\substack{i,j,k \geq 0\\
i+j+2k=6}} G_{ijk}x_0^ix_1^jy^k=0$,
where $q_x, q_y$ and $G_{ijk}$ are homogeneous polynomial on $\PP^1$.
\end{teo}

The proof of Theorem~\ref{main} requires some preparation, and we
postpone it to the end of the next section.

Let then $S$ be a minimal surface of general type having a free
canonical pencil of genus $2$ curves $f\colon S \rightarrow B$,
and let $X$ be the relative canonical model of $f$.

Let $\h$ be the horizontal fixed part of $|K_X|$, the union of
the components of $\Fix(|K_X|)$ which are not contracted by $f$.
Since, on a general fibre, $\h$ cuts a canonical divisor, the pull
back of a point of the canonical image, $\h=\gamma^*(\s)$ for some
section $\s$ of the conic bundle $\pi \colon C \rightarrow B$, to
which we can apply the results of the previous section.

$$
\xymatrix{%
\h\ar@{->>}[d]\ar@{^{(}->}[r]&X\ar@{->>}^{\gamma}[d]\ar@{->>}[r]^f&B\ar@{=}[d]\\
\s\ar@{^{(}->}[r] &C\ar@{->>}^{\pi}[r]&B\ar@/^1pc/[ll]\\
}
$$

\

\

Let $\K \subset \A$ be the ideal sheaf of the elements
vanishing along $\s$, $\SSS:=\A/\K$, and consider
the short exact sequence of vector bundles over $B$
\begin{equation}\label{k1}
0
\rightarrow
\K_1 \otimes \omega_B
\rightarrow
\A_1\otimes \omega_B=\R_1\otimes \omega_B =f_*\omega_S
\rightarrow
\SSS_1 \otimes \omega_B
\rightarrow
0.
\end{equation}

\begin{prop}\label{gangq}
Let $S$ be a minimal surface of general type having a free
canonical pencil of genus $2$ curves $f\colon S \rightarrow B$.
Let $b$ be the genus of $B$.

Then $0\leq \deg \SSS_1 \leq 1-q$. In particular $b \leq q(S) \leq 1$.
\end{prop}
\begin{Proof}
By a theorem of Fujita (\cite{fuj1}, \cite{fuj2}, \cite{zuchor},
see \cite[Remark 2.10]{cp})
$\R_1 \cong \hol_B^{q-b} \oplus A \oplus T$
where $A$ is ample and $T$ is a sum of torsion line bundles; in particular
$\deg(T)=h^0(T)=0$.

$S$ is of general type, and then $p_g \geq q \geq b$.
Assume by contradiction $p_g=b$. Then
$\deg \R_1=\chi -b+1=p_g-q-b+2=2-p_g \leq 0$ contradicts Fujita's theorem.
Therefore $p_g >b$.

By assumption every $2-$form on $S$ vanishes along $\h$,
and therefore the map $H^0(f_*\omega_S)
\rightarrow H^0(\SSS_1 \otimes \omega_B)$ induced by (\ref{k1})
is the zero map. Equivalently $h^0(\K_1\otimes \omega_B)=p_g$.
If $\K_1 \otimes \omega_B$ is nonspecial, $d:=\deg (\K_1 \otimes \omega_B) =p_g+b-1$. Else,
by Clifford's Theorem $d\geq 2(p_g-1)$.
By $p_g>b$ follows in both cases $d \geq  p_g+b-1$
and therefore
$$\deg \SSS_1=\deg \R_1-\deg \K_1 \leq  \chi -b+1-(p_g+b-1-2(b-1))=1-q$$

By the above mentioned theorem of Fujita, $\deg \SSS_1\geq 0$ and therefore $q\leq 1$.

\qed
\end{Proof}

By Proposition~\ref{s->tau'},
$\SSS$ is determined up to isomorphism by $\SSS_1$ and
a divisor $\tau' < \tau$.

\begin{lem}\label{comp}
$\tau'\neq 0$ and $\s \not\subset \Delta$.
\end{lem}
\begin{Proof}
Assume by contradiction $\s \cap P =\emptyset$. Then, in a
neighbourhood of $\h$, $\hol_X(K_X)$ is the pull back of
$\hol_{C\setminus P}(K_C+\zeta)$ where $\hol_{C\setminus P}(2\zeta)\cong
\hol_{C\setminus P}(\Delta)$. Therefore $\h K_X=2 \s (K_C+\zeta)$.

Since $K_{X|\PP^1}$ is the relative $\hol(1)$ of
$\Proj(\R)$, $|K_X|$ is given by maps
$\omega_{B}^{-1} \rightarrow \R_1$, $|K_C+\zeta|$ by maps
$\omega_{B}^{-1} \rightarrow \A_1$, and $(K_C+\zeta)_{|\s}$ by
maps $\omega_{B}^{-1} \rightarrow \SSS_1$.

So $\s (K_C+\zeta)= \deg (K_C+\zeta)_{|\s} = \deg \SSS_1 - \deg
\omega_{B}^{-1}\leq 1-b-2(b-1)=b-1$.
It follows that $K_X \h\leq 2(b-1)\leq 0$.
By assumption $K_X$ nef, and therefore $b=1$ and $K_X\h=0$.
But then $\h$ is a rational curve dominating
the elliptic curve $B$, a contradiction.

Therefore $\s \cap P \neq \emptyset$, which by (\ref{tau'not0})
implies $\tau'\neq 0$, and, since $\s$ is irreducible, by
(\ref{PDelta}) follows $\s \not\subset \Delta$.

\qed \end{Proof}

Hence follows (cf. \cite[Corollary~1 of Theorem~5.1]{xiao})

\begin{teo}[Xiao's inequality] \label{gang}
Let $S$ be a minimal surface of general type having a free
canonical pencil of genus $2$ curves $f\colon S \rightarrow B$.
Then
\begin{align*}
K_S^2 &\geq 4\chi+6q-10+3\deg (\tau-\tau')=4p_g+2q-6+3\deg (\tau-\tau')\\
&\geq 4\chi+6q-10=4p_g+2q-6.
\end{align*}

In particular, if $K^2 \leq 4\chi-8$, then $b=q=0$ and $\tau=\tau'$.
\end{teo}
\begin{Proof}
Since $\s \not\subset \Delta$, then the composition of maps $(\R_3')^2 \stackrel{\delta}{\rightarrow} \A_6 \rightarrow \SSS_6$ is nonzero.
Therefore $\Hom((\R_3')^2,\SSS_6) \neq 0$ and

\begin{align*}
6(1-q)+3\deg \tau'&\geq 6\deg \SSS_1+3\deg \tau' = \deg \SSS_6 \\\
&\geq \deg (\R_3')^2=2(\deg \tau +\chi+1)
\end{align*}
and the inequality follows by the formula (see
\cite[Theorem~3]{horP}, \cite[Theorem~4.13]{cp}) $K_S^2=2\chi+\deg \tau -6$.

\qed \end{Proof}

\begin{remark}
Some of the ideas of this proof are already present in the original proof of Xiao. Indeed
the geometric idea is similar to the idea of Xiao's
proof; we have shown that $\deg \SSS_6 - \deg (\R_3')^2 =
\Delta \s$ is nonnegative. 

In the next section we will see how, considering the non reduced divisor
$2\s$, this approach gives formulae which cannot be obtained by
intersection arguments.
\end{remark}

\begin{remark}\label{polizzi}
Maybe $\tau=\tau'$ holds more generally, and not only when $K^2 \leq 4\chi-8$.

We write $K_X$ as $\Phi+|M|$, where $\Phi$ is a divisor and $|M|$ is a movable linear system. Since we are assuming that $|K_X|$ is composed with a free pencil, $|M|$ has no base points.
As pointed out to the author by Polizzi, each isolated fixed point of the involution
on $X$ is a base point of the canonical system of the fibre through
it, and therefore $P \in \Phi$.

Anyway, the condition  $\tau=\tau'$ is slightly stronger; if we assume
$\tau$ reduced (for the sake of simplicity) then $\tau=\tau'
\Leftrightarrow P \subset \s$, and in general $\s\subsetneq\gamma(\Phi)$.
\end{remark}

{\bf Assumption:} From now on we assume $K^2 \leq 4\chi-8$.

By Theorem \ref{gang} it follows that $\tau=\tau'$ and moreover that the base curve is $\PP^1$, which is the first part of Theorem~\ref{main}. Our next goal is the computation of $\R$.

\begin{lem}\label{R1}
$\A_1=\R_1 \cong \hol_{\PP^1} (1) \oplus \hol_{\PP^1}(p_g+1)$
\end{lem}

\begin{Proof}
$S$ has genus $2$ pencil and we are assuming $K_S$ composed with
$f$. {\it i.e.} $|K|=\Phi+|f^*L|$ so $f_* \hol_S(K_S)=f_*\hol_S(\Phi)
\otimes \hol_{\PP^1}(L)$. By $p_g=h^0(K)=h^0(L)$, $f_* \hol_S(K_S)
\cong \hol_{\PP^1}(a) \oplus \hol_{\PP^1}(p_g-1)$ with $a<0$.

Finally $q(S)=0$ implies $h^1(f_* \hol_S(K_S))=0$ and therefore $a=-1$.

\qed \end{Proof}

We write
\begin{equation}\label{A1}
\A_1 \cong \hol_{\PP^1}(1)x_0 \oplus \hol_{\PP^1}(p_g+1)x_1
\end{equation}
labelling for our convenience its summands with the variables
$x_0,x_1$.
By Proposition~\ref{gangq}, $\SSS_1 \cong \hol_{\PP^1}(1)$ and
$x_1$ belongs to $\K_1$, the kernel of the map $\A_1 \rightarrow \SSS_1$. In fact
\begin{lem}\label{S2}
$\SSS \cong \A/x_1$.
\end{lem}
\begin{Proof}
Since $x_1$ belongs to the kernel of the map $\A \rightarrow \SSS$ we
have a surjection $\A/x_1 \rightarrow \SSS$, and we conclude by
showing that each graded piece of $\A/x_1$ is a line bundle.

Looking at the stalks of $\A$ as described in Remark~\ref{tau} this
fails if and only if there is a point $p\in\supp \tau$ such that $t$
divides $f_2(x_0,0;t)$ (here $f_2$ is the one in (\ref{stalk})). It
follows that the corresponding parameter $r''$ (introduced in
Remark~\ref{tau'}) is nonzero and then $\tau \neq \tau'$,
contradicting Corollary~\ref{gang}.

\qed \end{Proof}

In particular, denoting by $q\colon \A_2 \rightarrow \SSS_2$ the second graded piece of the surjection $\A \rightarrow \SSS$, we have an exact sequence
\begin{equation}\label{splits?}
0
\rightarrow
\hol_{\PP^1}(p_g+2)x_0x_1
\oplus
\hol_{\PP^1}(2p_g+2)x_1^2
\stackrel{(\sigma_2)_{|(x_1)}}{\longrightarrow}
\A_2
\stackrel{q}{\rightarrow}
\SSS_2
\rightarrow
0
\end{equation}
where $\SSS_2 \cong \hol_{\PP^1}(2+\deg \tau)=\hol_{\PP^1} (2p_g+\Theta)$ with
 $\Theta:=K^2-4p_g+6\in  \{0,1,2\}$.

$\A_2$ is a rank $3$ vector bundle over $\PP^1$, so there are $d_0\leq d_1
 \leq d_2$ with
\begin{equation*}
\A_2= \hol_{\PP^1}(d_0)y_0 \oplus  \hol_{\PP^1}(d_1)y_1 \oplus \hol_{\PP^1}(d_2)y_2
\end{equation*}

\begin{lem}\label{x12=y2}
$d_2=2p_g+2$. In particular, we can assume $y_2=\sigma_2(x_1^2)$.
\end{lem}
\begin{Proof}
The injectivity of $\sigma_2$ forces $d_2 \geq 2p_g+2$. If $d_2 >
2p_g+2$, then $d_2 > 2p_g + \Theta$ so $\hol_{\PP^1}(d_2)y_2 \subset
\ker q$. But by (\ref{splits?}) $\ker q \cong \hol_{\PP^1}(2p_g+2)
\oplus \hol_{\PP^1}(p_g+2)$, a contradiction.

\qed \end{Proof}

\begin{prop}\label{d0d1d2}
There exists $\alpha \geq 0$ such that
$$\begin{aligned}
d_0&=p_g+2+\alpha\\
d_1&=2p_g+\Theta - \alpha\\
d_2&=2p_g+2
\end{aligned}$$
\end{prop}
\begin{Proof}
We have computed $d_2$ in Lemma~\ref{x12=y2}; since by (\ref{A1}) and Lemma~\ref{A} $\deg \A_2=5p_g+\Theta+4$, we only have to prove $\alpha \geq 0$.

Else $d_1 > 2p_g+\Theta$, so $y_1 \in \ker q$. To embed
$\hol_{\PP^1}(d_1)y_1 \oplus \hol_{\PP^1}(2p_g+2)y_2 $ in $\ker q \cong
\hol_{\PP^1}(2p_g+2)x_1^2 \oplus \hol_{\PP^1}(p_g+2)x_0x_1$ we need $d_1 \leq
p_g+2$, so $2p_g+\Theta < p_g+2$, a contradiction to $\Theta \geq 0$,
$p_g \geq 2$.

\qed \end{Proof}

We consider the ordered bases $\{x_0^2,x_0x_1,x_1^2\}$ and $\{y_0,y_1,y_2\}$ of the source and of the target of $\sigma_2$. By Lemma~\ref{x12=y2}, $\sigma_2(x_1^2)=y_2$. Consider the homogeneous polynomials  $f_i$ and $g_j$ on $\PP^1$ such that $\sigma_2(x_0x_1)=\sum f_iy_i$, $\sigma_2(x_0^2)=\sum g_iy_i$.
The exact sequence (\ref{splits?}) gives
\begin{equation}\label{ficoprime}
\gcd(f_0,f_1)=1.
\end{equation}
Since $K_S^2 \leq 4\chi-8$, then $\Theta \leq 2$ and $\deg f_2=p_g \geq p_g+\Theta-2 =\deg f_0+ \deg f_1$. From (\ref{ficoprime}) follows that $f_2$ is in the ideal generated by $f_0$ and $f_1$. We may then assume, up to a suitable change of the coordinates in the target, that $f_2=0$. The matrix of $\sigma_2$ in the chosen bases is
\begin{equation}\label{sigma2}
\sigma_2=\begin{pmatrix}
g_0&f_0&0\\
g_1&f_1&0\\
g_2&0&1\\
\end{pmatrix}.
\end{equation}

\begin{prop}\label{alpha=0}
If $\alpha=0$, then $\R \cong (\Sym V)/\I$, with
$$
V:=\hol_{\PP^1}(1)x_0 \oplus \hol_{\PP^1}(\chi)x_1 \oplus \hol_{\PP^1}(K^2-2\chi+8)y \oplus \hol_{\PP^1}(K^2-\chi+7)z,
$$
the grading of $\Sym V$ is given by $\deg x_i=1,\ \deg y=2,\ \deg
z=3$, and $\I$ is the sheaf of ideals generated by
$x_0^2+q_x x_1^2+q_y y, z^2 +\sum_{\substack{i,j,k \geq 0\\
i+j+2k=6}} G_{ijk}x_0^ix_1^jy^k=0$,
where $q_x, q_y$ and $G_{ijk}$ are homogeneous polynomial on $\PP^1$.
\end{prop}
\begin{Proof}
If $\alpha=0$, then $\deg f_0=d_0-(p_g+2)=0$ and by
(\ref{ficoprime}) we can assume $f_0=1$.
Then the exact sequence (\ref{splits?}) splits and
therefore by Lemma~\ref{A}, $\A$ is the quotient of $\Sym \W$, where
$$\W \cong \hol_{\PP^1}(1)x_0 \oplus \hol_{\PP^1}(p_g+1)x_1 \oplus \hol_{\PP^1}(2p_g+\Theta)
y_1$$
by the ideal generated by
\begin{equation}\label{x0^2}
x_0^2-(g_0x_0x_1+(g_1-f_1g_0)y_1+g_2x_1^2).
\end{equation}
By changing the splitting of $\A_1$ in (\ref{A1}) , {\it i.e.} by changing the choice of $x_0$, we can assume $g_0=0$.

The statement follows, setting $y:=y_1$, from Proposition~\ref{Delta}.
Indeed the term $x_0^2$ in the equation (\ref{x0^2}) guarantees that
each map $\Sym^d \W \rightarrow \A_d$ splits.
In particular the map $\delta$ lifts to $\Sym^6 \W$
giving the equation $\sum G_{ijk} x_0^ix_1^jy^k$.

\qed \end{Proof}

\section{The sheaf of algebras of \protect $2\s$ and a lifting lemma}\label{2s}

\begin{df}
We consider the sheaf of algebras $\SSS':=\A/x_1^2$.
Since $\SSS=\A/x_1$, the
surjection $\A \rightarrow \SSS'$ correspond to the inclusion of the
non reduced divisor $2\s \subset C$.
\end{df}

$\A_6$ is a quotient of $\Sym^3 \A_2$, the cokernel of (see \cite[Lemma~4.4]{cp})
$$
i_3 \colon (\det \A_1)^2 \otimes \A_2 \rightarrow \Sym^3 \A_2
$$
defined by $$
i_3((x_0\wedge x_1)^2 \otimes v)=
(\sigma_2(x_0^2)\sigma_2(x_1^2)-\sigma_2(x_0x_1)^2)v.
$$
$\SSS_6'$ is the quotient of $\A_6$ by the multiples of $x_1^2$.
Killing first $y_2=x_1^2$ and then the multiples of
$\sigma_2(x_0x_1)^2=(f_0y_0+f_1y_1)^2$, we obtain $\SSS_6'$ as
cokernel of the map

$\begin{array}{cccc}
&&&\hol_{\PP^1}(3d_0)y_0^3\\
&&&\oplus\\
\hol_{\PP^1}(2p_g+4+d_0)&&&\hol_{\PP^1}(2d_0+d_1)y_0^2y_1\\
\oplus&\rightarrow&\Sym^3 \SSS_2' \cong&\oplus\\
\hol_{\PP^1}(2p_g+4+d_1)&&&\hol_{\PP^1}(d_0+2d_1)y_0y_1^2\\
&&&\oplus\\
&&&\hol_{\PP^1}(3d_1)y_1^3
\end{array}$

given by the matrix
$$
\begin{pmatrix}
f_0^2&0\\
2f_0f_1&f_0^2\\
f_1^2&2f_0f_1\\
0&f_1^2\\
\end{pmatrix}
$$
and therefore by (\ref{ficoprime}) $\SSS_6'$ is locally free of rank
$2$ and more precisely
$$
\SSS_6' \cong \hol_{\PP^1}(5p_g+2\Theta+\alpha+2)
\oplus \hol_{\PP^1}(6p_g+3\Theta-\alpha)
$$
and
the surjection $\Sym^3\SSS_2' \rightarrow \SSS_6'$ is given by the
matrix

\begin{equation}\label{S3S2->S6}
\begin{pmatrix}
3f_1^2&-2f_0f_1&f_0^2&0\\
0&f_1^2&-2f_0f_1&3f_0^2
\end{pmatrix}
\end{equation}

It follows
\begin{prop}\label{alphatheta}
$0\leq \alpha \leq \Theta$.
\end{prop}

\begin{Proof}
The surjection $\A \rightarrow \SSS$ factors through
$\SSS'$. By Lemma~\ref{comp} $\Hom(\R_3'^2,\SSS_6)\neq 0$ which implies
$\Hom(\R_3'^2,\SSS_6')\neq 0$, so
$2(3p_g+\Theta)\leq 6p_g+3\Theta-\alpha \Leftrightarrow \alpha \leq \Theta$.

\qed \end{Proof}

\begin{cor}\label{Theta0}
If $\Theta=0$, then $\R$ is the sheaf described in Proposition~\ref{alpha=0}.
\end{cor}
\begin{Proof}
If $\Theta=0$, by Proposition~\ref{alphatheta} follows $\alpha=0$ and we
can apply Proposition~\ref{alpha=0}.

\qed \end{Proof}

If $\Theta>0$, $\alpha$ may be different from $0$. The main point of the proof of Proposition~\ref{alpha=0} is that if $\alpha=0$ we can lift $\delta$ from $\A_6$ to a simpler target ($\Sym^6 \W$). When $\alpha>0$ we can't do that in general.

We look for a line bundle $\M$  and a map $0\neq m\in \Hom(\R_3'^2 \otimes \M^{-1},\R_3'^2)$
such that $\delta \circ m$ lifts to $\Sym^3 \A_2$, giving a commutative diagram

\begin{equation}\label{lifting}
\xymatrix{
\R_3'^2 \otimes \M^{-1}\ar^{\bar{\delta}}[r]\ar_{m}[d]&\ar^{}[d]\Sym^3 \A_2\\
\R_3'^2\ar^{\delta}[r]&\A_6\\
}
\end{equation}

Applying the functor $\Hom(\R_3'^2\otimes \M^{-1},\cdot)$ to the exact sequence
$$
0
\rightarrow
(\det \A_1)^2 \otimes \A_2
\stackrel{i_3}{\rightarrow}
\Sym^3 \A_2
\rightarrow
\A_6
\rightarrow
0
$$
we see that if $\Ext^1(\R_3'^2\otimes \M^{-1},(\det \A_1)^2 \otimes \A_2)=0$ ({\it i.e.} if $\deg \M$ is big enough) then the lift $\overline{\delta}$ exists.

Once we have such an $m$, to compute $\Hom(\R_3'^2,\A_6)$ we need to study
 which $\bar{\delta} \in \Hom
 (\R_3'^2 \otimes \M^{-1},\Sym^3 \A_2)$ belong to a diagram as
 (\ref{lifting}).
To simplify this computation, we need $\M$ of degree as small as
 possible. The forthcoming Lemma~\ref{liftlem} shows that we can take
 $m=f_0^4$, so $\deg \M=4\alpha$.

\begin{df}
Consider the natural decomposition $$\Sym^3 \A_2 = \sum_{i+j+k=3} \LLL_{ijk} y_0^{i}y_1^{j}y_2^{k}$$ as sum of vector bundles.
Clearly $\LLL_{ijk}\cong \hol(id_0+jd_1+kd_2)$. Then
$$
\V_{i\leq 1}:= \bigoplus_{i\leq 1} \LLL_{ijk} y_0^{i}y_1^{j}y_2^{k} \subset \Sym^3 \A_2.
$$
Similarly we define $\V_{i\geq 2}$; clearly $\Sym^3 \A_2=\V_{i\leq 1} \oplus \V_{i\geq 2}$.
\end{df}

\begin{lem}[Lifting Lemma]\label{liftlem}
Consider a map $\LLL \rightarrow \A_6$ where $\LLL$ is any line bundle. Then the composition map
$\LLL(-4\alpha) \stackrel{\cdot f_0^4}{\rightarrow} \LLL \rightarrow \A_6$
lifts to $\V_{i\leq 1} \subset \Sym^3 \A_2$.
\end{lem}
\begin{Proof}

Let $\J \subset \Sym \A_2$ be the kernel of the map $\Sym \A_2 \rightarrow  \A_{even} \subset \A$. The stalk of $\J$ at each point is a principal ideal generated by  the class of
\begin{equation}\label{Qfinal}
\Q:=(f_0y_0+f_1y_1)^2-y_2\sum g_iy_i.
\end{equation}

Therefore the kernel $\J_3$ of the map $\Sym^3\A_2 \rightarrow \A_6$
is generated by $\Q y_0, \Q y_1, \Q y_2$. By (\ref{ficoprime}), $f_0
\neq 0$, and then there are no nonzero relative cubics in the kernel without the terms $y_0^3,y_0^2y_1,y_0^2y_2$: $\V_{i\leq 1} \cap \J_3 = \{0\}$.
Defining $\T:=\coker \left( \J_3 \oplus \V_{i\leq 1} \rightarrow \Sym^3 \A_2 \right)$ we get the following commutative diagram of exact sequences
\begin{equation}\label{cubo}
\xymatrix{
&0\ar_{}[d]&0\ar_{}[d]&&\\
0\ar_{}[r]&\J_3\ar_{}[r]\ar_{}[d]&\J_3 \oplus \V_{i\leq 1}\ar_{}[r]\ar_{}[d]&\V_{i\leq 1}\ar_{}[r]\ar@/_1pc/[l]&0\\
&\Sym^3 \A_2\ar_{}[d]\ar_{}@{=}[r]&\V_{i\geq 2}\oplus \V_{i\leq 1}\ar_{}[d]&&\\
&\A_6\ar_{}[r]\ar_{}[d]&\T\ar_{}[r]\ar_{}[d]&0&\\
&0&0&&\\
}
\end{equation}
from which we deduce that $\V_{i\leq 1}$ maps isomorphically to $\ker\left( \A_6 \rightarrow \T\right)$.
In other words a map to $\A_6$ can be lifted to $\V_{i\leq 1}$ if and only its image goes to $0$ on $\T$: we conclude if we prove that $\T$ is annihilated by $f_0^4$.

Simplifying the term $\V_{i\leq 1}$ in the second column of
(\ref{cubo}), we see $\T$ as cokernel of a map between $\J_3$ and the
$\V_{i\geq 2}$, both sum of $3$ line bundles.
Since
$$
\begin{aligned}
\Q y_0&=f_0^2y_0^3+2f_0f_1y_0^2y_1-g_0y_0^2y_2+\ldots\\
\Q y_1&=f_0^2y_0^2y_1+\ldots\\
\Q y_2&=f_0^2y_0^2y_2+\ldots
\end{aligned}
$$
(here we have explicited the terms divisible by $y_0^2$, which are the
only terms relevant for this computation)
$$
\T \cong \coker \begin{pmatrix}
f_0^2&0&0\\
2f_0f_1&f_0^2&0\\
-g_0&0&f_0^2
\end{pmatrix}.
$$
Therefore $f_0^4$ annihilates $\T$.

\qed \end{Proof}

The map $\Sym \A_2 \rightarrow \A$ induces an inclusion $C \subset \Proj(\Sym
\A_2)$  and more precisely $C=\{(f_0y_0+f_1y_1)^2=y_2\sum g_iy_i\}$.
By Lemma~\ref{liftlem}, $\Delta \subset C$ is cut by
\begin{equation}\label{Grat}
F_{\Delta}:=\frac{\sum_{i\leq 1}F_{ijk}y_0^iy_1^jy_2^k}{f_0^4}
\end{equation}
where $F_{ijk}$ are polynomials on $\PP^1$ of degree
$$\deg F_{ijk} =(id_0+jd_1+kd_2)-(6p_g+2\Theta)+4\alpha,$$ and
moreover, since $F_{\Delta}$ has no poles when restricted to $C$
\begin{equation}\label{Gratmodf0^4}
\sum_{i\leq 1}F_{ijk}y_0^iy_1^jy_2^k=\left(\sum B_iy_i\right)\Q\
\text{ mod } f_0^4
\end{equation}
where the $B_i$ are rational functions on $\PP^1$ whose denominators
are invertible modulo $f_0$, and $\Q$ is as in (\ref{Qfinal}).

\begin{lem}\label{G120}
If $p_g>\alpha+2$ then $F_{120}=0$ and $F_{030}$ is nonzero and
divisible by $f_0^2$. Moreover $\gcd(f_0,g_0)=1$.
\end{lem}
\begin{Proof}
We compute the nonzero map $(\R_3')^2 \rightarrow \SSS_6'$ by
composing the map $\bar{\delta}\colon (\R_3')^2(-4\alpha)
\rightarrow \Sym^3 \A_2$ given by the numerator of (\ref{Grat}) with
the restriction $\Sym^3 \A_2 \rightarrow \Sym^3 \SSS_2'$, and then
with the matrix (\ref{S3S2->S6}), and finally dividing the result by
$f_0^4$.

We find that the map $(\R_3')^2 \rightarrow \SSS_6'$ is represented by
the matrix
$(\frac{F_{120}}{f_0^2},3\frac{F_{030}}{f_0^2}-2\frac{f_1F_{120}}{f_0^3})
$.

It follows that $f_0^2|F_{120}$. The assumption $p_g>\alpha+2$ is equivalent
to $\deg F_{120} < 2\alpha$, so $F_{120}=0$ and $F_{030}$ must be
nonzero and divisible by $f_0^2$.

Assume that there is a point $p$ with $g_0(p)=f_0(p)=0$. Then (see
(\ref{sigma2})) $p \in \supp \tau$. The node of $C$ above $p$, say
$P_0$, belongs to $\s=\{y_2=f_0y_0+f_1y_1=0\}$ and therefore has
relative coordinates  $(y_0:y_1:y_2)=(1:0:0)$.

Near $P_0$, $f_1 \neq 0$ by (\ref{ficoprime}), so to restrict $F_\Delta$ to $\s$ in a neighbourhood of $P_0$ we can substitute $y_2=0$ and
$y_1=-y_0f_0/f_1$: an equation for $\Delta_{|\s}$ is
$$
\frac{(-f_0^3F_{030}+f_0^2f_1F_{120})y_0^3}{f_0^4f_1^3}=-\frac{F_{030}}{f_0}\frac{y_0^3}{f_1^3}.
$$
Since $f_0^2|F_{030}$, $P_0\subset \Delta$, contradicting (\ref{PDelta}).

\qed \end{Proof}

Now we study the $B_i$ in (\ref{Gratmodf0^4}).
\begin{lem}\label{B}
If $p_g>\alpha+2$, then, modulo $f_0^4$,
\begin{align*}
B_0&=-2f_0\frac{F_{030}}{f_1^3}\\
B_1&=\frac{F_{030}}{f_1^2}\\
f_0^2B_2&=-2f_0\frac{g_0F_{030}}{f_1^3}
\end{align*}
Moreover $F_{111}$ is nonzero and divisible by $f_0^2$.
\end{lem}
\begin{Proof}
The expressions of the $B_i$ follow by looking at the
coefficients of (\ref{Gratmodf0^4}): $B_1$ by the term $y_1^3$, $B_0$
by $y_0y_1^2$ (since $F_{120}=0$ by Lemma~\ref{G120}), and $f_0^2B_2$
by $y_0^2y_2$.

Finally, looking at the coefficient of $y_0y_1y_2$, we find that
modulo $f_0^4$
\begin{align*}
f_0F_{111}=&-f_0(B_0g_1+B_1g_0)+2B_2f_0^2f_1\\
&=f_0\frac{F_{030}}{f_1^3} (2f_0g_1-f_1g_0-4f_1g_0)\\
&=uf_0F_{030}
\end{align*}
where by Lemma~\ref{G120}, $\gcd(u,f_0)=1$.

Since $f_0^2$ divides $F_{030}$, then it divides $F_{111}$ too. Moreover, since
$0\neq F_{030}=hf_0^2$ with $h$ of degree $\Theta-\alpha \in \{0,1\}$, then
$f_0F_{030}$ cannot vanish modulo $f_0^4$ and therefore $F_{111}$ is nonzero
too.

\qed \end{Proof}

As a consequence we find the following
\begin{cor}\label{alpha>0}
If $\alpha >0$ then $p_g \leq 2\alpha - \Theta+4$.
\end{cor}
\begin{Proof}
We can assume $p_g > \alpha+2 $ and apply Lemma~\ref{B}
to find $2\alpha=\deg f_0^2 \leq \deg F_{111}$ which is equivalent to
$p_g \leq 2\alpha - \Theta+4$.

\qed \end{Proof}

{\it Proof of Theorem~\ref{main}}.
Xiao's inequality guarantees that $S$ is regular
and the pencil has rational base.
If $\Theta=0$, Corollary~\ref{Theta0} ensures that $\R \cong (\Sym
V)/\I$. Else, by Proposition~\ref{alphatheta} and
Corollary~\ref{alpha>0} $p_g > \Theta+4 \geq 2\alpha - \Theta+4$ implies
$\alpha =0$ and we conclude by Proposition~\ref{alpha=0}.

\qed

The inequality in Corollary~\ref{alpha>0} is sharp. Indeed its proof
shows a method for constructing more examples with $p_g=2\alpha-\Theta+4$
by fixing $F_{111}/f_0^2$ and consequently computing $F_{030}$ and the
$B_i$.

\begin{example}\label{example}
We give examples for each pair $(\alpha,\Theta)$ with $1\leq \alpha \leq \Theta \leq 2$, $p_g=2\alpha-\Theta+4$: three examples with $(\alpha, \Theta, K^2,p_g)$ which equals respectively $(1,1,15,5)$, $(1,2,12,4)$ and $(2,2,20,6)$.

$f_0$ is a homogeneous polynomial on $\PP^1$ of degree $\alpha$: we choose
$$f_0:=
\begin{cases}
t_0 & \text{ if } \alpha=1\\
t_0(t_0-2) & \text{ if } \alpha=2\\
\end{cases}
$$
Note that $\deg g_1=2p_g+\Theta-\alpha-2=p_g+\alpha+2$ is even when
$\alpha=2$, so we can choose $g_1$ pure power of $f_0$.
We take
$$
f_1=t_1^{\alpha+2},\ g_0=t_1^{p_g+\alpha},\ g_1=f_0^{\frac{p_g+2}{\alpha}+1},\ g_2=0$$
so that $\tau=\{f_0^{\frac{p_g+2}{\alpha}+2}=t_1^{p_g+2\alpha+2}\}$.

The reader can check that
$$f_0(t_1^{\Theta-\alpha}y_1^2-y_0y_2)(f_0 y_1-2t_1^{p_g-2}y_2)-4f_0^2 y_0y_1y_2 $$
is, modulo $f_0^4$, of the form $(\sum B_iy_i)\Q$  (as prescribed by
(\ref{Gratmodf0^4}), note $F_{111}=-5f_0^2$).
So a candidate for the equation of $\Delta$ is
$$F_{\Delta_0}:=\frac{(t_1^{\Theta-\alpha}y_1^2-y_0y_2)(f_0 y_1-2t_1^{p_g-2}y_2)-4f_0 y_0y_1y_2}{f_0^3}.$$

The divisor $\Delta_0$ of the restriction of $F_{\Delta_0}$ to $C$ is
effective, but could be too
singular. Let $\D$ be the linear system of the effective
divisors on $C$ defined by the restriction of $F_{\Delta_0}+\lambda
y_2^3$, where $\lambda$ is any homogeneous polynomial on $\PP^1$ of
degree $6-2\Theta>0$. Since $\Fix \D \subset \{y_2=0\}
\cap C = 2\s$, its base points are contained
in $\s$. On the other hand the restriction of $F_{\Delta_0}$ to $\s$, computed
as in the proof of Lemma~\ref{G120}, is $t_1^{\Theta-\alpha}f_0$,
which has $\Theta$ distinct roots. So our pencil has $\Theta$ simple
base points, all along $\s$, and the general $\Delta \in \D$
is smooth by Bertini's theorem.
\end{example}

Finally we can give the proofs of Theorem~\ref{g=2} and of Theorem~\ref{allgenera?}. 

{\it Proof of Theorem~\ref{g=2}}.
$S$ is birational to the relative canonical model $X=\Proj(\R)$  of the
pencil. We have computed $\R$ in Theorem~\ref{main}. The description
as complete intersection in $\PP$ follows trivially.

Assume that $X$ is not the canonical model of $S$. Then there is an irreducible curve $E \subset X$ 
with $K_XE=0$. Since $X$ is the relative canonical model of $f$, then $E$ is not 
contracted by $f$.
Since $|K_X|=\Phi+|f^*L|$, then $E \cdot f^*L >0$ and therefore $E\Phi<0$, 
so $E \subset \Phi$ and more precisely $E\subset \h$.
For every fibre $F$ of $f$, $\h F=2$. Then either $E=\h$ or $\h=E+E'$ and $E$ and $E'$ are exchanged by the involution. 
In both cases $K_S\h=0$. But by the computations in Section~\ref{families}, 
$\h=\{x_1=0\}$ from which follows $0=K_S\h=(H-(p_g+1)F) \cdot (H-2F) \cdot 
Q \cdot G=2p_g+\Theta-4$, {\it i.e.} $\Theta=0$, $K^2=p_g=2$. 

All the other statements have been proved in Section~\ref{families}.
\qed

We need the following result of Sun.
\begin{teo}[\cite{sun1}, \cite{sun2}]\label{sun}
If $S$ is a minimal surface of general type whose canonical map is composed with a pencil of curves of genus $g$ without base points then $$\begin{matrix}
g=3&\Rightarrow&12 K^2 \geq 63 p_g-142\\
g=4&\Rightarrow&7K^2 \geq 48 p_g-134\\
g=5&\Rightarrow&9K^2 \geq 80 p_g-262\\
\end{matrix}
$$
\end{teo}

{\it Proof of Theorem~\ref{allgenera?}}.

If $S$ is a minimal surface of general type with a canonical pencil,
$K^2_S\leq 4\chi(\hol_S)-8$ and $\chi(\hol_S) \gg 0$, by Beauville's
Theorem~\ref{beauville} the pencil is free, and by
Theorem~\ref{sun} the pencil is a pencil of genus $2$ curves.
\qed

\end{document}